%% file: 623.tex
\documentstyle{article}

\def\bar#1{\overline{#1}}

\def\bar#1{\overline{#1}}

\def\bK{{\bf K}}

\input jbmac

\title{On the classifiability of Cellular Automata}
\author{John T. Baldwin
\thanks{Partially supported by NSF grant
9308768.}
\\ Department of Mathematics, Statistics and  Computer
Science\\University of Illinois at Chicago
\and
Saharon Shelah\thanks{This is paper 623 in Shelah's bibliography.  Both
authors thank
 the Binational Science Foundation for partial
support of this research.  Shelah's work on this project began during a
visit to
the University of Wisconsin partially supported by NSF \#144-EF67.}\\
Department of Mathematics\\
Hebrew University of Jerusalem\\ University of Wisconsin, Madison}

\begin{document}
\maketitle
\def\subi{\leq_{i}}
\def\subm{\leq_{s}}
\def\bT{{\bf T}}
\def\bK{{\bf K}}
\section{Background}
 Based on computer simulations Wolfram presented in several papers
conjectured classifications 
of cellular automata into 4 types.   
In \cite{Wolframct} Wolfram distinguishes the 4 classes of cellular
automata by the evolution
 of the pattern generated by applying a cellular automaton to a finite
input.  We quote from
page 161.
\begin{enumerate}
\item Pattern disappears with time.
\item Pattern evolves to a fixed finite size.
\item Pattern grows indefinitely at a fixed rate.
\item Pattern grows and contracts with time.
\end{enumerate}

Wolfram's qualitative classification is based on the examination of a large
number of 
simulations.  In addition to this classification based on the rate 
of growth, he
conjectured a similar classification according to the eventual pattern.  We
consider here
one formalization of his rate of growth suggestion. After completing our
major results (based only on Wolfram's work), 
we
investigated other contributions to the area and we report 
the relation of some them to our discoveries.  We thank Lyman Hurd, Nino
Boccaro, and Henryk Fuks for
their suggestions in this regard.

There are really two questions.  Can one classify the action of a cellular
automaton on a 
particular input $x$?  Can this be extended to a classification of automata
in terms, e.g., of
the average behavior of the automaton on all inputs?  It is straightforward
to prove  such a classification of
pairs $\langle A,x\rangle$.  That classification essentially is on the
lines Wolfram suggests.  (Only essentially, because class
three can be more precisely described as monotone growth.  The rate of
growth can vary from
$\log t$ to $t$.) But we show that this classification of pairs $(A,x)$
does not yield a classification of
automata $A$.  That is, for any nonnegative rationals $p,q$ with $p+q =1$,
we construct an automaton $A_{p,q}$ that depending on the input is 
likely to be in Class 3 with probability $p$ and in Class 4 with
probability $q$.  
In the process, we describe several patterns which seem to be qualitatively
different from those reported by Wolfram; in particular, with growth of order
 $\log t$.  
%Our work proceeds by encoding Turing machines as cellular
%automata in such a way that the space 
%used in the computation can be readily measured.  
We deal primarily with one dimensional cellular automata since they are
adequate for the counterexamples
we need.  The basic ideas extend naturally to higher dimensional cellular
automata.

There are a number of questions about the connections of these results with
Wolfram's conjectures.
First, Wolfram proposed several different schemes for classifying cellular
automata.  
Our result
argues that one of these conjectured classifications fails.  This does not,
a priori, invalidate
the other classifications.  
In particular, the formalization of the classification provided
by Culik and Yu \cite{CulikYu88} clearly divides all automata into 4
classes.  Similarly, that of \cite{Braga}
divides all automata into 3 classes.  
They show the question of which class a particular automata falls into is
undecidable.  Both of these formalizations classify an automaton by its
`worst case' or
most complicated behavior (as input varies).  Such a worst case
classification is completely consistent
with failure of an `average behavior' classification.     
Sutner \cite{Sutner} has shown the positions of the Culik-Yu classes in the
arithmetic hierarchy.  Like ours, these results deal with the action of
cellular automata on finite sequences; Ishii  \cite{Ishii} has established a
classification for the action on infinite sequences.

\section{Classification of automaton-input pairs}
\label{trueclass}
\jtbnot
\label{not1}  Our finite alphabets, usually denoted $\Sigma$, will always
contain the symbols 
$0,1,S,F,*,B$.  
$B$ will represent 
blank. A finite input on a two-way infinite tape will always be an initial
string of $B$'s, followed by
a finite word in $\Sigma$ (can include $B$'s) and then an infinite string
of $B$'s.
 We  assume that at
least one cell in an input string is not $B$.

\jtbdef  The {\em size} of a $1$-dimensional cellular automaton $A$ (or
Turing machine)
acting on input $x$ is the function,
$S_{A,x}(t)$ which assigns to each time $t$  the size of the configuration
on the tape
after $t$ steps of the computation with input $x$, that is the   distance
between 
left most and right most non-B cells at time $t$.

The following division into four cases is almost immediate from the
definitions.

\begin{lemm}   
\label{classAx} 
For every cellular automaton $A$ and finite input $x$, exactly one of the
following
holds.

\begin{enumerate}
\item $\lim_{t \rightarrow \infty} S_{A,x}(t) = 0$.

\item For some constant $c$, $0 < \limsup_{t \rightarrow \infty} S_{A,x}(t)
< c$.

\item $\lim_{t \rightarrow \infty} S_{A,x}(t) = \infty$ and  $S_{A,x}(t)$ is 
eventually monotone.

\item $\lim_ {t \rightarrow \infty} S_{A,x}(t) = \infty$ and  $S_{A,x}(t)$ is 
not eventually monotone.
\end{enumerate}
\end{lemm}

\Proof.  The key point is to note that if $\liminf_{t \rightarrow \infty}
S_{A,x}(t)$ is
bounded then a configuration is repeated and the action of $A$ on $x$ falls
into
the first or second class.$\Box$

We can refine this observation.

\begin{lemm}  Let $p = |\Sigma| +1$.  In either case 3 or 4 for a
$1$-dimensional
cellular automaton of radius $r$ we have:
 
$$\ln_p t \leq S_{A,x}(t) < |x| + rt.$$
\end{lemm}

\Proof.  The upper bound is immediate.  For, the lower bound note that if
$S_{A,x}(t)<\ln_p t$, two 
configurations must be repeated.  
But then the cellular automaton will cycle and we are in case 2. $\Box$

Note that the second class - pattern evolves to a fixed finite size -
encompasses both
periodic and {\em glider} configurations.  In a glider, the pattern repeats
cyclically
but moves across the domain.

\section{Simulation of Turing Machines by Cellular Automata}

We develop in this section a means first to simulate an arbitrary Turing
machine
on standard input by a cellular automata and then to simulate Turing
machines, 
which compute
total recursive functions, on arbitrary finite input.  
This gives rise to a convenient class of automata which we call
{\em dominating automata}.  Most discussions of recursive functions
are interested only in the computation of the value of a function
and input to the Turing machine is restricted to a standard form.  However,
Shepherdson \cite{Shepherdson65} dealt with arbitrary inputs
and proved, for example, that for every r.e. degree 
{\bf $\alpha$} there
is a Turing machine $M$ such that the collection of pairs of configurations
$\langle C_1, C_2\rangle$ such that $C_2$ appears if $M$ starts on $C_1$
has degree $\alpha$.

\jtbdef
\label{stdconfig}  
\begin{enumerate}
\item A {\em standard $\Sigma$-configuration} for a two-way tape contains a
unique
$S$ and $F$ and all non-$B$ cells are between them.
\item A {\em $\Sigma$-i/o (input/output)configuration} for a two-way tape
is a 
finite string surrounded by B's 
and beginning
with the symbol $S$ followed by a string of $0$'s and $1$'s 
(binary representation of a number m) 
followed by a *, 
followed by 
another string of $0$'s and $1$'s (binary representation of a 
number $n$) followed by an F.  
%An A may immediately precede the F.  
The strings may be empty.  
We write such a configuration
as $Sm*nF$.
\item We say a Turing machine is {\em on } a standard configuration if the 
head is reading a cell between the $S$ and $F$.
\end{enumerate}
\vskip .2in
A Turing machine is  specified by an alphabet $\Sigma$, a set  $Q$
of internal states, and a transition rule. The transition rule maps
the current state and the symbol  currently read to a new state,  
prints a symbol from $\Sigma$ and 
moves the head left or right.
It is easy to code the internal
states by expanding the alphabet; Turing machines thus become formally
more similar to cellular automata.  This coding is carried out more
precisely below.

\jtbnumpar{Operating Conventions}  We  restrict to Turing machines with
alphabet
$\Sigma$ which
obey the following.
\begin{enumerate}
\item  When reading $S$, the head can move left only if the $S$ is replaced by
a $0$ or $1$ (and S is printed on the cell to the left in
the next step).
\item  When reading $F$, the head can move right only if the $F$ is
replaced by
a $0$ or $1$  (and F is printed on the cell to the right in
the next step).
\item $S$ or $F$ is printed only in one of these two ways.
\item A $*$ is printed only if on the immediately preceding step a $*$ has
been overprinted
with a $0$ or $1$.
\end{enumerate}.

With these conventions it is easy to check the following lemma.  

\begin{lemm} If a Turing machine begins on a standard
$\Sigma$-configuration then 
every successive configuration is a standard $\Sigma$-configuration.
\end{lemm}

\jtbdef  A {\em $1$-dimensional cellular automaton} acts on a 2-way
infinite tape.  
Each cell contains a symbol from
 a finite alphabet $\Sigma$ (possibly B).  
The automaton has {\em radius }$r$ if
the value of a cell at time $t+1$ depends on the value at time $t$
of the cell and its $r$ predecessors
and $r$ successors.  

We deal primarily with 
radius 1 rules which determine
the next value of a cell depending only on the current value of 
the cell and its left and right
neighbors. We require that state $B$ is quiescent, any cellular automaton
takes 
an input 
which is all $B$ to $B$. Thus,  beginning on finite input our tape
will always contain only a finite number of non-$B$ cells.

\jtbnumpar{The Simulation Language}
%Suppose the Turing machine $T$ computes the function $f$.  
Let the Turing machine $T$ have the alphabet
$\Sigma = \{S,F, 0,*,1,B,\}$, operations $O = \{L,R\}$ 
and states $Q = \{q_0, \ldots q_k\}$.  
Thus $T$ is given by a function $T:Q \times 
\Sigma \mapsto \Sigma \times O \times Q$.
 We will contruct a 1-dimensional cellular
automaton
 with radius 1 that simulates $T$.
Let $Q_1 = Q \union \{B\}$.
%\union \{q^*_i: q_i \in Q\}$.
Let $\Sigma_1 = \Sigma \times Q_1\times \{L,R,H,B\}$.

Thus each member of $\Sigma_1$ codes a symbol of the original language, 
a state of the
original machine and the head position of the Turing machine. (In effect,
this creates a 3 track tape.) 
When confusion is unlikely to ensue, we will describe only the
 projection of the tape onto one coordinate.  
Thus we may say the tape reads $Sm*nF$ to mean
 the sequence of non-B first coordinates.
The cell is {\em active} if the head position (third coordinate) is H.
%The $q^*_i$ and $q'$ are auxilliary states which occur only in 
%nonstandard configurations and
%whose role  
%will be discussed in the proof of Theorem~\ref{findca}.
 We clarify this description
with the following definition.  Note that while $\Sigma$ and $\Sigma$-i/o 
configurations involved only the symbols from $\Sigma$, in
$\Sigma_1$-configurations
we also code the head position and state (of the simulated Turing machine). 

\jtbdef
\begin{enumerate}
\item
A {\em standard $\Sigma_1$-configuration} of a two-way infinite tape
satisfies the following.
\begin{enumerate}
\item All but finitely many cells  contain $B' = \langle B,B,B \rangle$. 
\item The first coordinates of the non-$B'$ cells form a standard
$\Sigma$-configuration.
\item There is a unique cell whose third coordinate is $H$;  all non-$B'$
cells to the
left of it have $R$ as third coordinate; all non-$B'$ cells to the
right of it  have $L$ as third coordinate.
\end{enumerate}
\item A standard $\Sigma_1$-configuration is a {\em standard $\Sigma_1$-i/o
configuration }
if in addition
\begin{enumerate}
\item  The first coordinates form a standard $\Sigma$-i/o configuration.
\item  One cell contains the entry $\langle S, q_o, H \rangle$.
The head position is $L$ for 
all other non-$B'$ cells.
\end{enumerate}
\end{enumerate}

Now we show that there is a simple computation of each partial recursive
function by a
cellular automaton.  This is, of course, well known.  For example, the
basic idea of the 
simulation here occurs in \cite{Hurd87}[6.3]. 
The argument we give here clarifies
and motivates our later constructions.  We see now how our simulation works
on standard
input.  Later, we introduce further complications to deal with nonstandard
input.

%ote that any cellular automata $A$ which simulated a universal Turing
% machine would have to%
%behave in each class defined in Lemma~\ref{classAx} for some inputs.  That
is, 
%a universal %automata
 %could not be charactized by a statement of the form, 
%for all inputs $x$, $\langle A, x\rangle$
%has a  rate of growth in, e.g. class 4.

The following argument is similar to the simulation described independently
but earlier in
\cite{CulikYu88}.  The novelty of the simulation in this paper 
appears in the treatment of nonstandard configurations.  
The important role of initial conditions and the possibility of 
 heads was pointed out independently but earlier
in \cite{Dharetal}; our analysis of this situation is new.

\begin{lemm}  
\label{sim1} For every Turing machine $M$, there is a cellular automaton
$A_M$ such that
the action of $M$,
beginning in state $q_0$ at the symbol $S$,  
on a standard $\Sigma$-i/o configuration is exactly the action of $A_M$ on
the 
first coordinates of the associated  standard $\Sigma_1$-i/o configuration.
Moreover, if the operation of $A_M$ begins on any standard
$\Sigma_1$-configuration,
all later configurations are also standard $\Sigma_1$-configurations.
\end{lemm}

\proof.  The automaton $A_M$ has dimension and radius 1. We describe the
action of $A_M$ on a 
cell
$i$ based on cells
$i-1$, $i$, $i+1$ (the site at $i$) with contents (for $ j = i-1, i, i+1$):
$\langle symbol, 
state, 
head \ position \rangle =
\langle s_j, q_j,p_j \rangle$.  The description here  is for cells which
appear in 
a standard $\Sigma_1$-configuration.  The definition is extended to
nonstandard 
configurations below. 
\begin{enumerate}
\item The first coordinate (i.e. the symbol)
at the next stage is determined entirely by cell $i$.
\begin{enumerate}
\item If $p_i \neq H$, the first coordinate remains the same.
\item If $p_i = H$, the first coordinate becomes the symbol printed by $M$
in state $q_i$
reading $s_i$.
\end{enumerate}
\item If $p_{i}$ is  $H$, the new state and head position of cell $i$ is
determined by 
cell $i$.  The state remains the same; the head position is $L$ or $R$
depending 
on whether $M$ moves left or right when reading $s_i$ in state $q_i$.
\item If $p_{i+1}$ is  $H$, the new state and head position of cell $i$ is
determined by 
cell $i+1$.
\begin{enumerate}
\item If in state $q_{i+1}$
reading $s_{i+1}$, $M$ moves left and goes into state $q'$, the new
position of cell $i$
is $H$ and the new state is $q'$. (The new position of cell $i+1$ is $L$.)
\item If in state $q_{i+1}$
reading $s_{i+1}$, $M$ moves right, the  position of cell $i$
remains $R$ and the state remains the same.  (The new position of cell
$i+1$ is $R$.)
\end{enumerate}
\item If $p_{i+1}$ is  $R$, then the new head position is again $R$ and the
state and symbol
are also unchanged.
\item If neither cell $i$, nor $i+1$ has head position $H$ or $R$, the new
state and 
head position of cell $i$ depend on 
cell $i-1$.
\begin{enumerate}
\item If $p_{i-1}$ is  $L$, then the new head position is $L$ and the state
and symbol
are unchanged.
\item If $p_{i-1} = H$, the new state and head position of cell $i$ is
determined by 
cell $i-1$.
\begin{enumerate}
\item If in state $q_{i-1}$
reading $s_{i-1}$, $M$ moves right and goes into state $q'$, the new
position of cell $i$
is $H$ and the new state is $q'$.  (The new position of $i-1$ is $R$.)
\item If in state $q_{i-1}$
reading $s_{i-1}$, $M$ moves left, the new position of cell $i$
is $L$ and the state remains the same.  (The new position of $i-1$ is $L$.)
\end{enumerate}

\end{enumerate}
\end{enumerate}
Just checking, one sees that on a standard $\Sigma_1$ i/o configuration, 
the simulation works as desired.$\Box$
\vskip .2in

We want to deal with arbitrary inputs.  We will arrange that on a finite
input, the
rightmost active cell will eventually dominate the computation.  In order
to do this
we have to restrict to certain kinds of computations of total recursive
functions.

\jtbnumpar{Normal input-output conventions}  The Turing machine $M$ {\em
normally computes}
 the function $f$, if for each $m$, on 
input $SmF$, beginning on $S$ in initial state $q_0$, it computes $Sf(m)F$
and halts.

We want to consider a 
 nonstandard input-output convention.

\jtbdef 
\begin{enumerate}
\item The Turing machine $M$ is said to {\em copy/compute} $f$ if beginning
at $S$ 
on a tape with
standard configuration $Sm*F$, it computes $Sm*f(m)F$ and halts.
\item
The Turing machine $T$ is said to {\em fully compute}
 the function $f$ on empty imput
%if placed with the head at any point in a standard 
%configuration $Sm*kF$,
 if the machine successively computes the sequences $Sn*f(n)F$
 for each natural number $n$.
%greater than $m$. 
%The symbol $A$ indicates the computation of $f$  at a particular 
%value has been completed.
\end{enumerate}

Obviously, every total recursive $f$ can be copy/computed by a 
Turing machine $T_f$.  The next remark is equally obvious; we spell
it out because we make use of the details in our simulation.

%\jtbdef  A Turing machine $T$ is {\em careful} if on every input $x$,
%it prints an $S$ somewhere and  returns infinitely often to each 
%cell containing $S$.

\begin{lemm}  
\label{carefultm} For any total recursive function $f$, there is a 
 Turing machine
$T_f$ which  fully computes $f$ on empty input.
\end{lemm}

\proof.  Fix a Turing machine $M$ which  
normally computes $f$.  
Now we describe the operation of the
new machine $T_f$ which fully computes $f$ on empty input.
  We assume the initial state is $q_o$. Using
special states it writes $S0*F$.  Now we begin the main loop.
It moves left erasing as it
goes until it reaches *.  
It then moves left adding 1 to the number  on the left of * and
moving the $S$ one cell to the left if necessary.  
The configuration now begins $Sm+1*$.
 Then head moves right and copies $m+1$ after *.  
Now it behaves on the sequence $*m+1$ as 
$M$ behaves on $Sm+1$ to compute $f(m+1)$.  
When it reaches the halting state of $M$,
this finishes one iteration of the loop. $\Box$

Note that since in incrementing $m$, $S$ is pushed to left (every $2^n$
steps) 
and $m$ is copied
to the right pushing $F$ to the right 
(unless the computation is already longer than $\log m$), any finite interval 
containing the initial configuration 
will eventually lie between $S$ and $F$.

We need one more refinement on our Turing computations; 
its use in this context was suggested to us by
Gyorgy Turan.  By an {\em initial position} of a Turing machine, 
we mean an input string,
a position of the head on the tape, and an initial state.

\begin{lemm}
\label{nocycle}
For any Turing machine $T$, there is a Turing machine $T'$ in a language
$\Sigma'$,
which on standard input simulates $T$, but does not cycle on any initial
position.
Moreover, on any $\Sigma'$-input $x$, 
$S_{T',x}(t) = \max (S_{T',x|\Sigma}(t),t)$.  Moreover for $t$ bigger
than the length of the input, $S_{T',x}(t)$ is a strictly increasing function.
\end{lemm}

\Proof.  Let $\Sigma'$ add a second track to the tape. The only symbols 
which occur on the
second track are  $0$,$1$ and  $B$. In accordance with the convention in 
Notation~\ref{not1}, this track contains only a finite 
number of $1$'s and $0$'s.
 $T'$ acts as $T$ on the first track.  At each step in the computation
the machine prints a $0$ on the 2nd track at the position currently being
read.  If that cell was blank it replaces the $0$ with $1$ and
proceeds to the next step of the computation. Otherwise, it moves to the
right until it reaches the first cell not $1$, prints a $1$ on it,
returns to the $0$, changes it to $1$ and 
proceeds to the next step of the computation. 
%k, it moves to the right reading
%the second 
%Between each step in the $\Sigma$-computation, adds one to the  
%string of $0$ and $1$'s (viewed as a binary number) on the
%second track that includes the head position and is bounded by the closed
blanks
%before and after the head position.
$\Box$

\jtbdef The cellular automaton $A$ is said to {\em completely compute}
 the total function $f$, if
for some $m$, the machine successively computes the sequences 
$Sn*f(n)F$ for each natural number $n>m$.

%\begin{enumerate}
%\item beginning in a standard $\Sigma_1$-configuration with first
coordinates $Sm*xF$,
% the machine successively computes the sequences $Sn*f(n)F$ for each
natural number $n>m$.
%\item beginning in a nonstandard configuration,
% the machine successively computes the sequences $Sn*f(n)F$ for each natural
%number $n$.
%\end{enumerate}

The initial input may contain a correct partial computation of $f(m)$;
in this case the machine just continues the computation.

\begin{thm}
\label{findca}
  For any total recursive function $f$, there is a cellular automaton
$A_f$ which  completely computes $f$.
\end{thm}

{\bf Proof Outline.}  Fix $T_f$, a Turing machine, which fully computes
$f$, and
which, using Lemma~\ref{nocycle}, does not cycle on any input.

We will establish two properties of the action of the simulating automaton.
\begin{enumerate}  
\item The successor of a standard $\Sigma_1$-configuration  $C$ is a
standard $\Sigma_1$-configuration $C'$.  Moreover, the first coordinates of 
$C'$ are the result of the action of $T_f$ on the first coordinates of $C$.
\item Any tape input with only finitely many non-$B'$-cells will evolve in
finitely many 
steps to a standard $\Sigma_1$-configuration.
\end{enumerate}
Together these two facts yield the theorem.

The first property follows directly from Lemma~\ref{sim1}.    
For the second, 
 we regard a 
cellular automaton with alphabet $\Sigma_1$ as a number of heads 
each performing
 a $\Sigma$-computation.  We will arrange that the rightmost of these 
heads eventually
dominates the computation and computes $f$.
We would like to construct an automaton that acted independently of 
input and just
started completely computing $f$. But, we have to allow for the 
possibility that
the initial position is in the midst of a correct computation.  
An arbitrary configuration
may contain many heads; it may contain none.

%To deal with these difficulties we introduce the auxiliary states
%$q_i^*$.  
Consider first that 
the initial configuration contains only one non-blank cell which 
contains $\langle S, q_0,H\rangle$.  From such a site the machine proceeds 
to
fully compute $f$ as in Lemma~\ref{carefultm}.
% but writing the symbol $q_i^*$ when the original
%machine would write $q_i$.  
It will print $*$ once and
this $*$ will never move.  We call this the {\em generating subroutine}.
We must explain what happens when there are other nonblank cells.

We say a subsequence of a configuration (in particular a site)
 is {\em acceptable} if it occurs in a simulation (as in Lemma~\ref{sim1})
of a computation  
beginning on
the  standard $\Sigma_1$-i/o configuration associated with $S*F$. 
(If the middle cell of a site is $\langle S, q_0,H\rangle$
and the  pair of the
second two cells occur in such a simulation then the site is 
acceptable.)
A {\em stop} cell is
one
of  
$\langle S, q_0,R\rangle$  or $\langle S, q_0,H\rangle$.
% q^*_0,H\rangle$ or $\langle S, q^*_0,R\rangle$.
A site is  {\em quiet} if the right most cell is a stop cell;
the center cell becomes $\langle B,q,R\rangle$ where $q$ was the
current state.
Any other
site is called a {\em generating site} and the new entry of cell $i$ is 
$\langle S, q_0,H\rangle$. 
The operation of the machine on a cell which contains a head depends
on whether the site centered on the cell is acceptable, quiet, or generating.
If it is acceptable, the simulation continues as in Lemma~\ref{sim1}; the 
other actions have just been described.

Now we give a global picture of the operation of the automata.

\begin{enumerate}
\item
\relax From any generating site the machine begins the generating subroutine.
  This operation has
priority (writing over any other input) unless the head finds a stop cell 
to its right. 

 \item
If a site is acceptable and contains a head, this head will either trace
out a complete
computation of $f$ or find a stop cell to its right.
\end{enumerate}

In either  case when the computation finds a stop cell the left 
$H$ becomes $R$ and 
remains quiescent   until it
is eventually overwritten by the head on the right. 
(If there is a head on the right this will happen because the * written
by the right Head will never move; eventually the rightmost Head will write
over anything
written by the other Heads.)

To see that this machine computes a final sequence of values for $f$,
we analyze the initial string from the right.  Either the entire
configuration is acceptable or there is a right most generating site
followed by an acceptable string.  In the first case, the configuration
is a standard $\Sigma_1$-i/o configuration and the result follows by
Lemma~\ref{sim1}.  In the second case a complete computation of $f$ 
will propagate from the rightmost generating site. 
%using the $q^*_i$ states.  
The input to the right of this state will be used; the
input to the left is irrelevant to the eventual computation.$\Box$

\jtbdef We call an automaton $A_f$ constructed as in the proof of
Theorem~\ref{findca},
a {\em dominating automaton}.

Note that a dominating automaton uses unbounded space on any input, so  
a classification of automata according to the schema suggested
would have to put each dominating automaton in class 3 or 4.

\section{Composition and Nonclassifiability of Cellular Automata}
\label{classautin}

In this section we show how to compose a finite set of dominating cellular
automata 
$A_1 \ldots A_n$ into a single
automaton $A$ with a larger alphabet whose growth rates reflects that of
each $A_i$.  Moreover,
this composition can be chosen so that the classification of the behavior
of $A$ 
on input $x$ falls into specified type 3 or 4 with arbitrary probability. 

\jtbdef  Let $A_1, \ldots A_n$  be cellular automata of the same dimension
and radius with alphabet $\Sigma_0 \subseteq \Sigma$.  Let $X=
\bigcup_{i<n} X_i$ be an 
additional set of finite symbols (where the $X_i$ are disjoint). 
Form the language $\Sigma_1 = \Sigma \times X$.  Define the cellular automaton
$A = \oplus_i A_i$ with the following transition rule.  If the central cell
has an element
of $X_i$ as its second component use the transition rule from $A_i$ on the
first components.
$A$ is called the {\em composition of the $A_i$ with respect to $X$}.

We clearly have:

\begin{lemm}  If $A_1$ and $A_2$ are dominating automata then so is their
composition
(for any $X$).
\end{lemm}

\jtbdef  For each $n$, let $P_n$ be the probability measure assigning the
same probability
to each element of $\Sigma^n$  (i.e. each finite input of length $n$).

\jtbdef Let $P_n(i,A)$ be the probability that among all inputs $x$ of length
$n$, the function $S_{A,x}$ is in class i (from the classification in
Lemma~\ref{classAx}).
\vskip .2in

We now show that the  classification of Section~\ref{trueclass}
  does  not extend  from pairs $\langle
A, x \rangle$ to cellular automata $A$.

\begin{lemm}  
\label{noclass1}
Let $p,q$ be rational numbers $0 \leq p, q \leq 1$ with $p+q=1$.
There is a cellular automaton $A_{p,q}$ such that for every $n$,
$P_n(\{x: \langle A_{p,q}, x \rangle  \in {\rm Class\ 3}\}) = p$ and
$P_n(\{x: \langle A_{p,q}, x \rangle  \in {\rm Class \; 4}\}) = q$.
\end{lemm}

\Proof. Suppose there is such a classification. 
Let $A_1$ be in class 3 and $A_2$ in class $4$ represent two total
recursive functions
as in Theorem~\ref{findca}.
Choose $X$ with $p$ symbols for  $A_1$, $q$ symbols for $A_2$.  Now, the
required machine is the composition of the $A_i$ with respect to this
$X$.$\Box$

\jtbnumpar{Remark}  This construction  refutes  a  rigid classification
 of cellular automata into four classes according to the rate of growth
schema. 
It does not seem  to refute the separation of the bounded space
from unbounded space automata.  Two complications present themselves.  If the
automaton which is supposed to dominate is of class 1, it might die out before
it had a chance to exert its dominance over some pretenders.  This can be
remedied
by inserting a "resurrection state".  More seriously, if the "dominating
automaton"
were to glide to the right, it would never exert his dominance over, e.g. a
class
3 automaton to its left and we would be left with a class three pattern
instead 
of class 2.  This tends to support the judgement of \cite{Braga} who combine
classes 3 and 4 in their classification.
 
\section{Rate of Growth}

In  this section we investigate the rate of growth
of patterns generated by cellular automata.
The following examples shows that a pattern which grows monotonically in
size need not
grow at a `fixed rate' if that phrase is interpreted as `linearly in $t$'.

\jtbnumpar{A slow growing example}  Let $A$ be the cellular automaton which
is derived
 from the identity function by the construction in Lemma~\ref{findca}.
Then on standard
input $S*F$, $A$ successively writes $Sm*mF$ for any natural number $m$.
Thus, since it
takes time $\log n$ to write $n$,
$\lim_{t\rightarrow \infty}S_{A,S*F}(t)/\log t$ is a constant.

The difficulty of distinguishing the third and fourth classes is emphasized
by another
construction.

\jtbnumpar{Enforcing Monotonicity}
\label{monotuniv}  Let $A$ be any cellular automaton of class 3 or 4.  For
simplicity, suppose $\Sigma =\{0,1,B,S,F\}$ and that $A$ is $1$-dimensional
of radius 1.  
(The $S$, $F$ are inessential and included only
to keep our notation consistent.) We add a new symbol $M$ (for marked).  
Let $\Sigma'= \Sigma \times \{B,M\}$.
Let the value of $A'$ on three consecutive cells $\langle x_{i-1}, y_{i-1}
\rangle$,
 $\langle x_{i}, y_{i} \rangle$, $\langle x_{i+1}, y_{i+1} \rangle$
 be $\langle B,B\rangle$ if all the $x$'s and $y$'s are blank.  Otherwise
the second
coordinate is $M$ and the first coordinate is the result of applying $A$ to
$x_{i-1}, x_i, x_{i+1}$.   Then every cell that is ever
marked remains marked, so $A'$ is class 3 even if $A$ is class 4,
but the `information content' remains the same as that for $A$.

\jtbnumpar{Eventual Behavior}

The crux of the argument here is that the behavior of the function
$S_{A,x}(t)$ depends
essentially on both $A$ and $x$.  Paradoxically, we achieved this by
constructing automata
whose eventual behavior is independent of input in the following sense.

\begin{lemm} If $A_f$ is a cellular automaton from Lemma~\ref{findca} which
fully computes 
$f$.
\begin{enumerate}
\item For any input $x$, $\liminf _{t \rightarrow \infty} S_{A_f,x}(t) =
\infty$.
\item For any input $x$, there exist constants $t_0$ and $c$ such that for
$t \geq t_0$,
$$S_{A_f,x}(t) = S_{A_f,S*F}(t-c).$$
\end{enumerate}
Thus, the eventual behavior of $S_{A_f,x}(t)$ on any input $x$ is
determined by the eventual
behavior of $S_{A_f,S*F}$.
\end{lemm}

\proof.  If $x$ is nonstandard, after $t_0$ steps, $A_f$  settles on the
unique active cell, 
prints $S*F$, and simulates $T_f$ on input $S*F$.  
If $x$ is standard, the computation of $A_f$ on
$x$ begins $c$ steps into the computation of $A_f$ on $S$.
$\Box$

\jtbnumpar{Classifying minima}

Class 4 automata were defined by the property that $S_{A,x}(t)$ is 
not eventually monotone.  There are some restrictions on this nonmonotonicity.
For example,   for any $1$-dimensional cellular automaton, the
function which enumerates the points $(t_i,S_{A,x}(t_i))$ which are local
minima
%(local maxima) 
of $S_{A,x}(t)$ is clearly recursive.  We show that, in a certain sense, every
total
recursive function can be represented in this way.  Let $M$ be an
arbitrary Turing machine and $A_M$ be the cellular automata associated 
with $M$ in 
Lemma~\ref{sim1}.  Let $x$ be the input $\langle \langle S,q_0,H\rangle,
\langle S,0,L\rangle
\langle *,B,L\rangle, \langle F,B,L\rangle \rangle$ with all 
other cells $B'$.   Then for all $t$, 
the contents of the tape at time $t$ is the same whether considering
computation by $M$ or by $A_M$.  In particular, $S_{M,x}(t) =S_{A_M,x}(t)$.
For any total recursive function $f$, we
 constuct a Turing machine $M_f$ so that the contents of the tape at the
$2i$th 
minimum of $S_{M_f,x}(t)$ is  $Si*f(i)$.
We compute a total recursive function
$f$ by a Turing machine $M_f$ which uses strictly increasing space
 on 
the computation
of each value (as in Lemma~\ref{nocycle}).
 Note that 
space (using the second track) will strictly increase until the 
$\Sigma$-configuration
reads $Sn*f(n)$.  
When the computation is complete, add one more symbol to the second track.
Then erase the second track until it has the same length as the first.
(The interpolated step guarantees there is at least one step 
in this process.)  Again add one element to the second track, then erase
both tracks until the contents of the first are $Sn*F$.  Now, increment
$n$ to $n+1$ and compute $f(n+1)$; use the second track to guarantee that
the space used is increasing throughout this stage.  Thus the only 
space minina are at configuratons $Sn*F$ and $Sn*f(n)F$.  We have shown:

\begin{thm}  For any recursive function $f$ there is a Turing machine $M$ 
which 
computes $f$, and 
there is a cellular automata $A_M$ such that 
%if $m(i)$ denotes the ith local maximum
%of $S_{A_M,S0*F}(t)$, then 
%\begin{enumerate}
%\item 
the $\Sigma$-configuration of the $2i$th local minimum is $Si*f(i)F$.
%\item $m(i) - (\lg(i)+2)$ is the time taken by $M$ to compute $f(i)$.
%\end{enumerate}
\end{thm}

If we used 1-ary rather than binary notation we could easily decode
the value of $f$ directly from the values of $S_{A_M,S0*F}(t)$ at minima.

\section{Conclusions} We briefly compare these results with several 
related papers.

\jtbnumpar{Universality and Class}  Culik and Yu \cite{CulikYu88} gave a 
different formalization of
Wolfram's classification.  Paraphrasing slightly, they define

\begin{enumerate}
\item  $A$ evolves to all blanks from every finite input.
\item  $A$ has an ultimately periodic evolution on every finite input.
\item  For any two configurations $c_1$ and $c_2$, it is decidable
whether $c_1$ will evolve to $c_2$ under $A$.
\item  All other cellular automata.
\end{enumerate}

Clause four guarantees that this is a (cumulative) hierarchy classifying
all cellular automata.  
The spirit of this classification is to label each automaton with its most
complicated
behavior (ranging over all inputs).

Their Theorem 10 asserts that no universal automaton can be Class Three.
But our third class
is clearly a subset of theirs and we showed in Paragraph~\ref{monotuniv}
how to encode a
universal automaton into our third class.  The seeming paradox is resolved
by noting the
significance of input/output coding.  They report their result is obvious.
Indeed, it is
given that their i/o coding is unique.  That is, if (as specified in
\cite{CulikYu88}) there is a unique configuration representing each
natural number, then deciding whether $c_1$ evolves to $c_2$ under $A_f$ is
the same as
deciding whether $f$ on the input coded by $c_1$ gives the value coded by
$c_2$.  However,
in the scheme described in Paragraph~\ref{monotuniv} there are infinitely
many codes for each
possible output and so the contradiction is avoided. 

\jtbnumpar{Probabilities on infinite strings}  Ishii \cite{Ishii} has given
a probabalistic
classification of the behavior of cellular automata on infinite strings.
Informally,
an automata is in class X if for almost every intitial configuration (in a
specified measure
on $\Sigma^Z$) evolves to a configuration of type $X$.  While this result
is in a different
direction from ours, the distinction demonstrates again the importance of
distinguishing
behavior on finite strings from behavior on infinite strings.  An analogous
situation is the
contrast between
the undecidability of the ring of integers (arbitrary finite sequences) and
the decidability
of the field of real numbers (arbitrary sequences).

\jtbnumpar{The number of states}  In our construction, we freely expanded
the language $\Sigma$
by adding a small number of additional symbols.  The necessity of such an
expansion is made
clear by the proof by Land and Belew \cite{LandBelew} that
for any density $\rho$, there is no two-state automata (of any
radius) which can correctly decide whether sequences of arbitrary length
have density greater
than $\rho$.  In particular there can be no two state universal cellular
automata.  So
our use of more states was essential.

\jtbnumpar{Summary} This paper highlights the importance of input and
output conventions in describing the information content as opposed to the
dynamics
of a computation.  If  the automaton acts with the standard input/output
convention (\ref{stdconfig}), then a cellular
automaton simulating a universal Turing machine will, depending on the
input, have runs in each of the four classes.  However, by modifying the
output convention as in \ref{monotuniv}, we can construct a universal
cellular automaton which
behaves in class 3 on every input.   We have formalized Wolfram's
classification scheme in terms of the spatial rate of growth of a
computation.  We see that this notion is
well defined for pairs of an automaton
 acting on an input but that it can not be extended even probabalistically
to a classification of automata.  Several new patterns have been discovered
in the course of this investigation.  In one case the size increases
monotonically but
at a rate of $\log t$ rather than linearly.  
Wolfram describes class four automata as having complex localized structure
which is sometimes long lasting.   The examples of dominating class 4
given in this paper are different. 
 After a finite amount of chaotic (in a nontechnical sense) behavior 
they evolve to a pattern which grows monotonically on one side and as 
eratically as the time taken to compute  a given recursive function on the
other.

\bibliography{ssgroups}
\bibliographystyle{plain}
\end{document}

%% file: jbmac.tex
\makeatletter
\def\@begintheorem#1#2{\it \trivlist \item[\hskip
\labelsep{\bf #2\ #1.}]}
\def\@opargbegintheorem#1#2#3{\it \trivlist
      \item[\hskip \labelsep{\bf #1\ #2\ (#3).}]}
\makeatother

\newtheorem{them}{Theorem}[section]

\newtheorem{lemm}[them]{Lemma}

\newtheorem{thm}[them]{Theorem}

\newcommand{\jtbnumpar}[1]{\refstepcounter{them}
\trivlist
\item[\hskip \labelsep{\bf \thethm \ #1.}]}

\newcommand{\jtbdef}{\jtbnumpar{Definition}}
\def\jtbnot{\jtbnumpar{Notation}}
\newtheorem{ssnote}{COMMENT}
\def\PROOF #1.{\par\noindent{\it Proof#1}.\ \ignorespaces}
\def\proof #1.{\par\noindent{\it Proof#1}.\ \ignorespaces}

\def\subm{\leq}

\let\?=\joinrel

\let\leftv=^
\let\rightv=^

\def\k/{\kern.2em}    % for eg Theorem 1.2\k/ii)
\def\max{\mathop{\rm max}}

\let\union=\cup             %

\edef\bigcup{\mathop{\textstyle\mathchar\the\bigcup}}
\edef\bigcap{\mathop{\textstyle\mathchar\the\bigcap}}

\edef\bigwedge{\mathop{\textstyle\mathchar\the\bigwedge}}

\edef\bigvee{\mathop{\textstyle\mathchar\the\bigvee}}

\edef\sum{\mathop{\textstyle\mathchar\the\sum}}
\def\math&{\ \& \ }

\def\force {\mathrel^\joinrel\rightarrow}
\def\force {\mathrel{\scriptstyle\mathrel^\joinrel\rightarrow}}
\def\forceq {\mathrel{\mathop{\force}\limits_{\textstyle\texsim}}}
\def\forceq{\mathrel^\joinrel
 \mathrel{\mathop{\rightarrow}\limits_{\smash{\textstyle\texsim}}}}
\def\forceq{\mathrel{\scriptstyle\mathrel^\joinrel
 \mathrel{\smash{\mathop{\rightarrow}\limits_{\smash{\raise
 2pt\hbox{$\scriptstyle\texsim$}}}}}}}

\let\exclaim=!                 %
\let\texsim=\sim         %
\let\conj\sim             %
\def\conjp #1 {\conj_{#1}}     %
\let\sim\simeq             %
\def\0bar{\bar 0}         %
\def\1bar{\bar 1}         %

\def\hbar{\overline h}

\def\Proof{Proof}

\def\text#1{\ifmmode\leavevmode\hbox{#1}\else
   \typeout{Warning: \string\text \space used outside math mode!}
   \begingroup\hbox{#1}\endgroup\fi}